\documentclass{amsart}

\newtheorem{theorem}{Theorem}[section]

\newtheorem{corollary}[theorem]{Corollary}

\theoremstyle{definition}

\theoremstyle{remark}

\numberwithin{equation}{section}

\begin{document}
\title[Identities of symmetry for Euler polynomials]{$\begin{array}{c}
         \text{Identities of symmetry for Euler polynomials}\\
         \text{arising from quotients of fermionic integrals invariant under $S_{3}$}\
       \end{array}$}

\author{dae san kim and Kyoung ho park} \thanks{This work was supported by National Foundation of Korea Grant funded by the Korean Government(2009-0072514).
}
\address{Department of Mathematics, Sogang University, Seoul 121-742, Korea}
\email{dskim@sogong.ac.kr}
\address{Department of Mathematics, Sogang University, Seoul 121-742, Korea}
\email{sagamath@yahoo.co.kr}
\thanks{}

\subjclass[]{}

\date{}

\dedicatory{}


\begin{abstract}
In this paper, we derive eight basic identities of symmetry in three
variables related to Euler polynomials and alternating power sums.
These and most of their corollaries are new, since there have been
results only about identities of symmetry in two variables. These
abundance of symmetries shed new light even on the existing
identities so as to yield some further interesting ones. The
derivations of identities are based on the $p$-adic integral
expression of the generating function for the Euler polynomials and
the quotient of integrals that can be expressed as the exponential
generating function for the alternating power sums.

\end{abstract}

\maketitle

\section{Introduction and preliminaries}
  Let $p$ be a fixed odd prime. Throughout this paper,
$\mathbb Z_{p},\mathbb Q_{p},\mathbb C_{p}$ will respectively denote
the ring of $p$-adic integers, the field of $p$-adic rational
numbers and the completion of the algebraic closure of $\mathbb
Q_{p}$. For a uniformly differentiable(also called continuously
differentiable) function $f:\mathbb Z_{p}\rightarrow \mathbb C_{p}$,
the $p$-adic fermionic integral of $f$ is defined by
\begin{align*}
\int_{\mathbb Z_{p}}f(z)d\mu_{-1}(z)=\lim_{N\rightarrow
\infty}\sum_{j=0}^{p^{N}-1}{f(j)(-1)^{j}}.
\end{align*}
Then it is easy to see that
\begin{equation}\label{a}
\int_{\mathbb Z_{p}}{f(z+1)d\mu_{-1}(z)}+\int_{\mathbb
Z_{p}}{f(z)d\mu_{-1}(z)}=2f(0).
\end{equation}
 Let $|~|_{p}$ be the normalized absolute value of $\mathbb C_{p}$, such that $|p|_{p}= \frac{1}{p}$, and let
\begin{equation}\label{b}
E=\{t\in \mathbb C_{p}||t|_{p}<p^{-\frac{1}{p-1}}\}.
\end{equation}
Then, for each fixed $t \in E $, the function $f(z)=e^{zt}$  is
analytic on $\mathbb Z_{p}$ and by applying (\ref{a}) to this $f$,
we get the $p$-adic integral expression of the generating function
for Euler numbers $E_{n}$:
\begin{equation}\label{c}
\int_{\mathbb
Z_{p}}{e^{zt}d\mu_{-1}(z)}=\frac{2}{e^{t}+1}=\sum_{n=0}^{\infty}E_{n}
\frac{t^{n}}{n!}~(t \in E).
\end{equation}
So we have the following $p$-adic integral expression of the
generating function for the Euler polynomials $E_{n}(x)$:
\begin{equation}\label{d}
\int_{\mathbb
Z_{p}}{e^{(x+z)t}d\mu_{-1}(z)}=\frac{2}{e^{t}+1}e^{xt}=\sum_{n=0}^{\infty}E_{n}(x)
\frac{t^{n}}{n!}~(t \in E,x \in \mathbb Z_{p}).
\end{equation}
  Let $T_{k}(n)$ denote the alternating $k$th power sum of the first $n+1$ nonnegative integers, namely
\begin{equation}\label{e}
T_{k}(n)=\sum_{i=0}^{n}(-1)^{i}i^{k}=(-1)^{0}0^{k}+(-1)^{1}1^{k}+(-1)^{2}2^{k}+
\cdots +(-1)^{n}n^{k}.
\end{equation}
In particular,
\begin{equation}\label{f}
T_{0}(n)=
\begin{cases}
\begin{split}
1,& \quad \text {if} ~n\equiv 0~(mod~2),\\
0,& \quad \text {if} ~n\equiv 1~(mod~2),
\end{split}
\end{cases}
T_{k}(0)=
\begin{cases}
\begin{split}
1,& \quad \text {for} ~k=0,\\
0,& \quad \text {for} ~k>0.
\end{split}
\end{cases}
\end{equation}
From (\ref{c}) and (\ref{e}), one easily derives the following
identities: for any odd positive integer $w$,
\begin{equation}\label{g}
\frac{\int_{\mathbb Z_{p}}e^{xt}d\mu_{-1}(x)}{\int_{\mathbb
Z_{p}}e^{wyt}d\mu_{-1}(y)}=
\sum_{i=0}^{w-1}(-1)^{i}e^{it}=\sum_{k=0}^{\infty}T_{k}(w-1)\frac{t^{k}}{k!}~(t\in
E).
\end{equation}
In what follows, we will always assume that the $p$-adic fermionic
integrals of the various exponential functions on $\mathbb Z_{p}$
are defined for $t\in E$ (cf. (\ref{b})), and therefore it will not
be mentioned.

 \cite{DR1}, \cite{H1}, and \cite{T1}-\cite{Y1} are some of the previous
works on identities of symmetry involving Bernoulli polynomials and
power sums. These results were generalized in \cite{D1} to obtain
identities of symmetry involving three variables in contrast to the
previous works involving just two variables.

  In this paper, we will produce 8 basic identities of symmetry in three
variables $w_{1},w_{2},w_{3}$ related to Euler polynomials and
alternating power sums(cf. (\ref{p1}), (\ref{q1}), (\ref{t1}),
(\ref{w1}), (\ref{a2}), (\ref{c2}), (\ref{e2}), (\ref{f2})). These
and most of their corollaries seem to be new, since there have been
results only about identities of symmetry in two variables in the
literature. These abundance of symmetries shed new light even on the
existing identities. For instance, it has been known that (\ref{h})
and (\ref{i}) are equal and (\ref{j}) and (\ref{k}) are so(cf. [4,
Theorem 5, 7]). In fact, (\ref{h})-(\ref{k}) are all equal, as they
can be derived from one and the same $p$-adic integral. Perhaps,
this was neglected to mention in \cite{T1}. Also, we have a bunch of
new identities in (\ref{l})-(\ref{o}). All of these were obtained as
corollaries(cf. Cor.9, 12, 15) to some of the basic identities by
specializing the variable $w_{3}$ as 1. Those would not be unearthed
if more symmetries had not been available.

  Let $w_{1},w_{2}$ be any odd positive integers. Then we have:
\begin{align}
\label{h}
&\sum_{k=0}^{n}\binom{n}{k}E_{k}(w_{1}y_{1})T_{n-k}(w_{2}-1)w_{1}^{n-k}w_{2}^{k}\\
\label{i}
=&\sum_{k=0}^{n}\binom{n}{k}E_{k}(w_{2}y_{1})T_{n-k}(w_{1}-1)w_{2}^{n-k}w_{1}^{k}\\
\label{j}
=&w_{1}^{n}\sum_{i=0}^{w_{1}-1}(-1)^{i}E_{n}(w_{2}y_{1}+\frac{w_{2}}{w_{1}}i)\\
\label{k}
=&w_{2}^{n}\sum_{i=0}^{w_{2}-1}(-1)^{i}E_{n}(w_{1}y_{1}+\frac{w_{1}}{w_{2}}i)\\
\label{l}
=&\sum_{k+l+m=n}^{}\binom{n}{k,l,m}E_{k}(y_{1})T_{l}(w_{1}-1)T_{m}(w_{2}-1)w_{1}^{k+m}w_{2}^{k+l}\\
\label{m}
=&w_{1}^{n}\sum_{k=0}^{n}\binom{n}{k}\sum_{i=0}^{w_{1}-1}(-1)^{i}E_{k}(y_{1}+\frac{i}{w_{1}})T_{n-k}(w_{2}-1)w_{2}^{k}\\
\label{n}
=&w_{2}^{n}\sum_{k=0}^{n}\binom{n}{k}\sum_{i=0}^{w_{2}-1}(-1)^{i}E_{k}(y_{1}+\frac{i}{w_{2}})T_{n-k}(w_{1}-1)w_{1}^{k}\\
\label{o}
=&(w_{1}w_{2})^{n}\sum_{i=0}^{w_{1}-1}\sum_{j=0}^{w_{2}-1}(-1)^{i+j}E_{n}(y_{1}+\frac{i}{w_{1}}+\frac{j}{w_{2}}).
\end{align}

  The derivations of identities will be based on the $p$-adic integral
expression of the generating function for the Euler polynomials in
(\ref{d}) and the quotient of integrals in (\ref{g}) that can be
expressed as the exponential generating function for the alternating
power sums. We indebted this idea to the paper \cite{T1}.

\section{Several types of quotients of fermionic integrals}

  Here we will introduce several types of quotients of $p$-adic fermionic
integrals on $\mathbb Z_{p}$ or $\mathbb Z_{p}^{3}$ from which some
interesting identities follow owing to the built-in symmetries in
$w_{1},w_{2},w_{3}$. In the following, $w_{1},w_{2},w_{3}$ are all
positive integers and all of the explicit expressions of integrals
in (\ref{q}), (\ref{s}), (\ref{t}), and (\ref{v}) are obtained from
the identity in (\ref{c}).
\\
\\
(a) Type $\Lambda_{23}^{i}$ (for $i=0,1,2,3$)
\begin{align}
\label{p} I(\Lambda_{23}^{i})&=\frac{\int_{\mathbb
Z_{p}^{3}}e^{(w_{2}w_{3}x_{1}+w_{1}w_{3}x_{2}+w_{1}w_{2}x_{3}+w_{1}w_{2}w_{3}(\sum_{j=1}^{3-i}y_{j}))t}d\mu_{-1}(x_{1})d\mu_{-1}(x_{2})d\mu_{-1}(x_{3})}
{(\int_{\mathbb Z_{p}}e^{w_{1}w_{2}w_{3}x_{4}t}d\mu_{-1}(x_{4}))^{i}}\\
\label{q}
&=\frac{2^{3-i}e^{w_{1}w_{2}w_{3}(\sum_{j=1}^{3-i}y_{j})t}(e^{w_{1}w_{2}w_{3}t}+1)^{i}}
{(e^{w_{2}w_{3}t}+1)(e^{w_{1}w_{3}t}+1)(e^{w_{1}w_{2}t}+1)};
\end{align}
\\
\\
(b) Type $\Lambda_{13}^{i}$ (for $i=0,1,2,3$)
\begin{align}
\label{r} I(\Lambda_{13}^{i})&=\frac{\int_{\mathbb
Z_{p}^{3}}e^{(w_{1}x_{1}+w_{2}x_{2}+w_{3}x_{3}+w_{1}w_{2}w_{3}(\sum_{j=1}^{3-i}y_{j}))t}d\mu_{-1}(x_{1})d\mu_{-1}(x_{2})d\mu_{-1}(x_{3})}
{(\int_{\mathbb Z_{p}}e^{w_{1}w_{2}w_{3}x_{4}t}d\mu_{-1}(x_{4}))^{i}}\\
\label{s}
&=\frac{2^{3-i}e^{w_{1}w_{2}w_{3}(\sum_{j=1}^{3-i}y_{j})t}(e^{w_{1}w_{2}w_{3}t}+1)^{i}}
{(e^{w_{1}t}+1)(e^{w_{2}t}+1)(e^{w_{3}t}+1)};
\end{align}
\\
\\
(c-0) Type $\Lambda_{12}^{0}$
\begin{align}
\label{s0} I(\Lambda_{12}^{0})&=\int_{\mathbb
Z_{p}^{3}}e^{(w_{1}x_{1}+w_{2}x_{2}+w_{3}x_{3}+w_{2}w_{3}y+w_{1}w_{3}y+w_{1}w_{2}y)t}d\mu_{-1}(x_{1})d\mu_{-1}(x_{2})d\mu_{-1}(x_{3})
\\
\label{t}
&=\frac{8e^{({w_{2}w_{3}+w_{1}w_{3}+w_{1}w_{2}})yt}}{(e^{w_{1}t}+1)(e^{w_{2}t}+1)(e^{w_{3}t}+1)};
\end{align}
\\
\\
(c-1) Type $\Lambda_{12}^{1}$
\begin{align}
\label{u} I(\Lambda_{12}^{1})&=\frac{\int_{\mathbb
Z_{p}^{3}}e^{(w_{1}x_{1}+w_{2}x_{2}+w_{3}x_{3})t}d\mu_{-1}(x_{1})d\mu_{-1}(x_{2})d\mu_{-1}(x_{3})}
{\int_{\mathbb
Z_{p}^{3}}e^{(w_{2}w_{3}z_{1}+w_{1}w_{3}z_{2}+w_{1}w_{2}z_{3})t}d\mu_{-1}(z_{1})d\mu_{-1}(z_{2})d\mu_{-1}(z_{3})}
\\
\label{v}
&=\frac{(e^{w_{2}w_{3}t}+1)(e^{w_{1}w_{3}t}+1)(e^{w_{1}w_{2}t}+1)}{(e^{w_{1}t}+1)(e^{w_{2}t}+1)(e^{w_{3}t}+1)}.
\end{align}

  All of the above $p$-adic integrals of various types are invariant under
all permutations of $w_{1},w_{2},w_{3}$ as one can see either from
$p$-adic integral representations in (\ref{p}), (\ref{r}),
(\ref{s0}), and (\ref{u}) or from their explicit evaluations in
(\ref{q}), (\ref{s}), (\ref{t}), and (\ref{v}).
\section{Identities for Euler polynomials}
  In the following $w_{1},w_{2},w_{3}$ are all odd positive integers
except for (a-0) and (c-0), where they are any positive
integers.\\
\\
(a-0) First, let's consider Type $\Lambda_{23}^{i}$, for each
$i=0,1,2,3.$ The following results can be easily obtained from
(\ref{d}) and (\ref{g}).
\begin{equation*}
\begin{split}
I(\Lambda_{23}^{0})=\int_{\mathbb
Z_{p}}e^{w_{2}w_{3}(x_{1}+w_{1}y_{1})t}d\mu_{-1}(x_{1})
&\int_{\mathbb Z_{p}}e^{w_{1}w_{3}(x_{2}+w_{2}y_{2})t}d\mu_{-1}(x_{2})\quad\qquad\qquad\qquad\\
\times&\int_{\mathbb Z_{p}}e^{w_{1}w_{2}(x_{3}+w_{3}y_{3})t}d\mu_{-1}(x_{3})\\
\end{split}
\end{equation*}
\begin{equation*}
~\qquad=(\sum_{k=0}^{\infty}\frac{E_{k}(w_{1}y_{1})}{k!}(w_{2}w_3{}t)^k)
(\sum_{l=0}^{\infty}\frac{E_{l}(w_{2}y_{2})}{l!}(w_{1}w_{3}t)^l)
(\sum_{m=0}^{\infty}\frac{E_{m}(w_{3}y_{3})}{m!}(w_{1}w_{2}t)^m)\\
\end{equation*}
\begin{equation}\label{v0}
=\sum_{n=0}^{\infty}(\sum_{k+l+m=n}\binom{n}{k,l,m}E_{k}(w_{1}y_{1})E_{l}(w_{2}y_{2})E_{m}(w_{3}y_{3})w_{1}^{l+m}w_{2}^{k+m}w_{3}^{k+l})\frac{t^{n}}{n!},~\quad\\
\end{equation}
where the inner sum is over all nonnegative integers $k,l,m$, with
$k+l+m=n$, and
\begin{equation}\label{w}
\binom{n}{k,l,m}=\frac{n!}{k!l!m!}.
\end{equation}
\\
\\
(a-1) Here we write $I(\Lambda_{23}^{1})$ in two different ways:
\\
\\
(1) $I(\Lambda_{23}^{1})$
\begin{equation}\label{x}
\begin{split}
=\int_{\mathbb
Z_{p}}e^{w_{2}w_{3}(x_{1}+w_{1}y_{1})t}d\mu_{-1}(x_{1})
&\int_{\mathbb Z_{p}}e^{w_{1}w_{3}(x_{2}+w_{2}y_{2})t}d\mu_{-1}(x_{2})\quad\qquad\qquad\qquad\\
\times&\frac{\int_{\mathbb
Z_{p}}e^{w_{1}w_{2}x_{3}t}d\mu_{-1}(x_{3})}{\int_{Z_{p}}e^{w_{1}w_{2}w_{3}x_{4}t}d\mu_{-1}(x_{4})}
\end{split}
\end{equation}
\begin{equation*}
\qquad=(\sum_{k=0}^{\infty}{E_{k}(w_{1}y_{1})\frac{(w_{2}w_{3}t)^k}{k!}})
(\sum_{l=0}^{\infty}{E_{l}(w_{2}y_{2})\frac{(w_{1}w_{3}t)^l}{l!}})
(\sum_{m=0}^{\infty}{T_{m}(w_{3}-1)\frac{(w_{1}w_{2}t)^m}{m!}})\quad\quad
\end{equation*}
\begin{equation}\label{y}\qquad=\sum_{n=0}^{\infty}(\sum_{k+l+m=n}\binom{n}{k,l,m}E_{k}(w_{1}y_{1})E_{l}(w_{2}y_{2})T_{m}(w_{3}-1)w_{1}^{l+m}w_{2}^{k+m}w_{3}^{k+l})\frac{t^{n}}{n!}.\qquad\quad
\end{equation}
\\
\\
(2) Invoking (\ref{g}), (\ref{x}) can also be written as
\\
\\
\begin{equation*}
I(\Lambda_{23}^{1})=\sum_{i=0}^{w_{3}-1}(-1)^i\int_{\mathbb
Z_{p}}e^{w_{2}w_{3}(x_{1}+w_{1}y_{1})t}d\mu_{-1}(x_{1})
\int_{\mathbb Z_{p}}e^{w_{1}w_{3}(x_{2}+w_{2}y_{2}+\frac{w_{2}}{w_{3}}i)t}d\mu_{-1}(x_{2})\\
\end{equation*}
\begin{equation*}
=\sum_{i=0}^{w_{3}-1}(-1)^i(\sum_{k=0}^{\infty}E_{k}(w_{1}y_{1})\frac{(w_{2}w_{3}t)^k}{k!})
(\sum_{l=0}^{\infty}E_{l}(w_{2}y_{2}+\frac{w_{2}}{w_{3}}i)\frac{(w_{1}w_{3}t)^l}{l!})\\
\end{equation*}
\begin{equation}\label{z}
=\sum_{n=0}^{\infty}(w_{3}^{n}\sum_{k=0}^{n}\binom{n}{k}E_{k}(w_{1}y_{1})\sum_{i=0}^{w_{3}-1}(-1)^iE_{n-k}(w_{2}y_{2}+\frac{w_{2}}{w_{3}}i)w_{1}^{n-k}w_{2}^{k})\frac{t^{n}}{n!}.\quad\qquad\\
\end{equation}
(a-2) Here we write $I(\Lambda_{23}^{2})$ in three different ways:
\\
\\
(1) $I(\Lambda_{23}^{2})$
\begin{equation}\label{a1}
=\int_{\mathbb
Z_{p}}e^{w_{2}w_{3}(x_{1}+w_{1}y_{1})t}d\mu_{-1}(x_{1}) \times
\frac{\int_{\mathbb
Z_{p}}e^{w_{1}w_{3}x_{2}t}d\mu_{-1}(x_{2})}{\int_{\mathbb
Z_{p}}e^{w_{1}w_{2}w_{3}x_{4}t}d\mu_{-1}(x_{4})} \times
\frac{\int_{\mathbb
Z_{p}}e^{w_{1}w_{2}x_{3}t}d\mu_{-1}(x_{3})}{\int_{\mathbb
Z_{p}}e^{w_{1}w_{2}w_{3}x_{4}t}d\mu_{-1}(x_{4})}
\end{equation}
\begin{equation*}
=(\sum_{k=0}^{\infty}{E_{k}(w_{1}y_{1})\frac{(w_{2}w_{3}t)^k}{k!}})
(\sum_{l=0}^{\infty}{T_{l}(w_{2}-1)\frac{(w_{1}w_{3}t)^l}{l!}})
(\sum_{m=0}^{\infty}{T_{m}(w_{3}-1)\frac{(w_{1}w_{2}t)^m}{m!}})\quad\quad
\end{equation*}
\begin{equation}\label{b1}
=\sum_{n=0}^{\infty}(\sum_{k+l+m=n}\binom{n}{k,l,m}E_{k}(w_{1}y_{1})T_{l}(w_{2}-1)T_{m}(w_{3}-1)w_{1}^{l+m}w_{2}^{k+m}w_{3}^{k+l})\frac{t^{n}}{n!}.\qquad\quad
\end{equation}
\\
\\
(2) Invoking (\ref{g}), (\ref{a1}) can also be written as
\\
\\
~$I(\Lambda_{23}^{2})$
\begin{equation}\label{c1}
=\sum_{i=0}^{w_{2}-1}(-1)^i\int_{\mathbb
Z_{p}}e^{w_{2}w_{3}(x_{1}+w_{1}y_{1}+\frac{w_{1}}{w_{2}}i)t}d\mu_{-1}(x_{1})
\times \frac{\int_{\mathbb
Z_{p}}e^{w_{1}w_{2}x_{3}t}d\mu_{-1}(x_{3})}{\int_{\mathbb
Z_{p}}e^{w_{1}w_{2}w_{3}x_{4}t}d\mu_{-1}(x_{4})}~
\end{equation}
\begin{equation*}
=\sum_{i=0}^{w_{2}-1}(-1)^i(\sum_{k=0}^{\infty}E_{k}(w_{1}y_{1}+\frac{w_{1}}{w_{2}}i)\frac{(w_{2}w_{3}t)^k}{k!})
(\sum_{l=0}^{\infty}T_{l}(w_{3}-1)\frac{(w_{1}w_{2}t)^l}{l!})
\end{equation*}
\begin{equation}\label{d1}
=\sum_{n=0}^{\infty}(w_{2}^{n}\sum_{k=0}^{n}\binom{n}{k}\sum_{i=0}^{w_{2}-1}(-1)^iE_{k}(w_{1}y_{1}+\frac{w_{1}}{w_{2}}i)T_{n-k}(w_{3}-1)w_{1}^{n-k}w_{3}^{k})\frac{t^{n}}{n!}.\quad
\end{equation}
\\
\\
(3) Invoking (\ref{g}) once again, (\ref{c1}) can be written as
\\
\\
~$I(\Lambda_{23}^{2})$
\begin{equation*}
=\sum_{i=0}^{w_{2}-1}\sum_{j=0}^{w_{3}-1}(-1)^{i+j}\int_{\mathbb
Z_{p}}e^{w_{2}w_{3}(x_{1}+w_{1}y_{1}+\frac{w_{1}}{w_{2}}i+\frac{w_{1}}{w_{3}}j)t}d\mu_{-1}(x_{1})~
\end{equation*}
\begin{equation*}
=\sum_{i=0}^{w_{2}-1}\sum_{j=0}^{w_{3}-1}(-1)^{i+j}\sum_{n=0}^{\infty}E_{n}(w_{1}y_{1}+\frac{w_{1}}{w_{2}}i+\frac{w_{1}}{w_{3}}j)\frac{(w_{2}w_{3}t)^n}{n!}
\end{equation*}
\begin{equation}\label{e1}
=\sum_{n=0}^{\infty}((w_{2}w_{3})^n\sum_{i=0}^{w_{2}-1}\sum_{j=0}^{w_{3}-1}(-1)^{i+j}E_{n}(w_{1}y_{1}+\frac{w_{1}}{w_{2}}i+\frac{w_{1}}{w_{3}}j))\frac{t^{n}}{n!}.\quad
\end{equation}
\\
\\
(a-3)
\\
\\
~$I(\Lambda_{23}^{3})$
\begin{equation*}
=\frac{\int_{\mathbb
Z_{p}}e^{w_{2}w_{3}x_{1}t}d\mu_{-1}(x_{1})}{\int_{\mathbb
Z_{p}}e^{w_{1}w_{2}w_{3}x_{4}t}d\mu_{-1}(x_{4})}
\times\frac{\int_{\mathbb
Z_{p}}e^{w_{1}w_{3}x_{2}t}d\mu_{-1}(x_{2})}{\int_{\mathbb
Z_{p}}e^{w_{1}w_{2}w_{3}x_{4}t}d\mu_{-1}(x_{4})}
\times\frac{\int_{\mathbb
Z_{p}}e^{w_{1}w_{2}x_{3}t}d\mu_{-1}(x_{3})}{\int_{\mathbb
Z_{p}}e^{w_{1}w_{2}w_{3}x_{4}t}d\mu_{-1}(x_{4})}
\end{equation*}
\begin{equation*}
=(\sum_{k=0}^{\infty}T_{k}(w_{1}-1)\frac{(w_{2}w_{3}t)^{k}}{k!})
(\sum_{l=0}^{\infty}T_{l}(w_{2}-1)\frac{(w_{1}w_{3}t)^{l}}{l!})
(\sum_{m=0}^{\infty}T_{m}(w_{3}-1)\frac{(w_{1}w_{2}t)^{m}}{m!})
\end{equation*}
\begin{equation}\label{f1}
=\sum_{n=0}^{\infty}\sum_{k+l+m=n}^{}(\binom{n}{k,l,m}T_{k}(w_{1}-1)T_{l}(w_{2}-1)T_{m}(w_{3}-1)w_{1}^{l+m}w_{2}^{k+m}w_{3}^{k+l})\frac{t^{n}}{n!}.\quad
\end{equation}
\\
\\
(b) For Type $\Lambda_{13}^{i}~(i=0,1,2,3)$, we may consider the
analogous things to the ones in (a-0), (a-1), (a-2), and (a-3).
However, these do not lead us to new identities. Indeed, if we
substitute $w_{2}w_{3},w_{1}w_{3},w_{1}w_{2}$ respectively for
$w_{1},w_{2},w_{3}$ in (\ref{p}), this amounts to replacing $t$ by
$w_{1}w_{2}w_{3}t$ in (\ref{r}). So, upon replacing
$w_{1},w_{2},w_{3}$ respectively by
$w_{2}w_{3},w_{1}w_{3},w_{1}w_{2}$ and dividing by
$(w_{1}w_{2}w_{3})^n$, in each of the expressions of
(\ref{v0}),(\ref{y}), (\ref{z}), (\ref{b1}), (\ref{d1})-(\ref{f1}),
we will get the corresponding symmetric identities for Type
$\Lambda_{13}^{i}~(i=0,1,2,3)$.
\\
\\
(c-0)
\\
\begin{align*}
&I(\Lambda_{12}^{0})\\
&=\int_{\mathbb Z_{p}}e^{w_{1}(x_{1}+w_{2}y)t}d\mu_{-1}(x_{1})
\int_{\mathbb Z_{p}}e^{w_{2}(x_{2}+w_{3}y)t}d\mu_{-1}(x_{2})\int_{\mathbb Z_{p}}e^{w_{3}(x_{3}+w_{1}y)t}d\mu_{-1}(x_{3})\\
&=(\sum_{n=0}^{\infty}\frac{E_{k}(w_{2}y)}{k!}(w_{1}t)^{k})(\sum_{l=0}^{\infty}\frac{E_{l}(w_{3}y)}{l!}(w_{2}t)^{l})(\sum_{m=0}^{\infty}\frac{E_{m}(w_{1}y)}{m!}(w_{3}t)^{m})
\end{align*}
\begin{align}\label{g1}
=\sum_{n=0}^{\infty}(\sum_{k+l+m=n}\binom{n}{k,l,m}E_{k}(w_{2}y)E_{l}(w_{3}y)E_{m}(w_{1}y)w_{1}^{k}w_{2}^{l}w_{3}^{m})\frac{t^n}{n!}.~\qquad\qquad
\end{align}
\\
\\
(c-1)
\\
\begin{equation*}
\frac{\int_{\mathbb
Z_{p}}e^{w_{1}x_{1}t}d\mu_{-1}(x_{1})}{\int_{\mathbb
Z_{p}}e^{w_{1}w_{2}z_{3}t}d\mu_{-1}(z_{3})}
\times\frac{\int_{\mathbb
Z_{p}}e^{w_{2}x_{2}t}d\mu_{-1}(x_{2})}{\int_{\mathbb
Z_{p}}e^{w_{2}w_{3}z_{1}t}d\mu_{-1}(z_{1})}
\times\frac{\int_{\mathbb
Z_{p}}e^{w_{3}x_{3}t}d\mu_{-1}(x_{3})}{\int_{\mathbb
Z_{p}}e^{w_{3}w_{1}z_{2}t}d\mu_{-1}(z_{2})}
\end{equation*}
\begin{equation*}
~\qquad=(\sum_{k=0}^{\infty}T_{k}(w_{2}-1)\frac{(w_{1}t)^{k}}{k!})
(\sum_{l=0}^{\infty}T_{l}(w_{3}-1)\frac{(w_{2}t)^{l}}{l!})
(\sum_{m=0}^{\infty}T_{m}(w_{1}-1)\frac{(w_{3}t)^{m}}{m!})
\end{equation*}
\begin{equation}\label{h1}
=\sum_{n=0}^{\infty}(\sum_{k+l+m=n}^{}\binom{n}{k,l,m}T_{k}(w_{2}-1)T_{l}(w_{3}-1)T_{m}(w_{1}-1)w_{1}^{k}w_{2}^{l}w_{3}^{m})\frac{t^{n}}{n!}.\quad
\end{equation}
\section{Main theorems}
  As we noted earlier in the last paragraph of Section 2, the various
types of quotients of $p$-adic fermionic integrals are invariant
under any permutation of $w_{1},w_{2},w_{3}$. So the corresponding
expressions in Section 3 are also invariant under any permutation of
$w_{1},w_{2},w_{3}$. Thus our results about identities of symmetry
will be immediate consequences of this observation.

  However, not all permutations of an expression in Section 3 yield
distinct ones. In fact, as these expressions are obtained by
permuting $w_{1},w_{2},w_{3}$ in a single one labelled by them, they
can be viewed as a group in a natural manner and hence it is
isomorphic to a quotient of $S_{3}$. In particular, the number of
possible distinct expressions are $1,2,3,$ or $6$. (a-0), (a-1(1)),
(a-1(2)), and (a-2(2)) give the full six identities of symmetry,
(a-2(1)) and (a-2(3)) yield three identities of symmetry, and (c-0)
and (c-1) give two identities of symmetry, while the expression in
(a-3) yields no identities of symmetry.

  Here we will just consider the cases of Theorems 8 and 17, leaving
the others as easy exercises for the reader. As for the case of
Theorem 8, in addition to (\ref{v1})-(\ref{x1}), we get the
following three ones:
\begin{align}\label{i1}
\sum_{k+l+m=n}\binom{n}{k,l,m}E_{k}(w_{1}y_{1})T_{l}(w_{3}-1)T_{m}(w_{2}-1)w_{1}^{l+m}w_{3}^{k+m}w_{2}^{k+l},
\end{align}
\begin{align}\label{j1}
\sum_{k+l+m=n}\binom{n}{k,l,m}E_{k}(w_{2}y_{1})T_{l}(w_{1}-1)T_{m}(w_{3}-1)w_{2}^{l+m}w_{1}^{k+m}w_{3}^{k+l},
\end{align}
\begin{align}\label{k1}
\sum_{k+l+m=n}\binom{n}{k,l,m}E_{k}(w_{3}y_{1})T_{l}(w_{2}-1)T_{m}(w_{1}-1)w_{3}^{l+m}w_{2}^{k+m}w_{1}^{k+l}.
\end{align}
But, by interchanging $l$ and $m$, we see that (\ref{i1}),
(\ref{j1}), and (\ref{k1}) are respectively equal to (\ref{v1}),
(\ref{w1}), and (\ref{x1}).
\\
As to Theorem 17, in addition to (\ref{f2}) and (\ref{g2}), we have:
\begin{align}
\label{l1}
&\sum_{k+l+m=n}\binom{n}{k,l,m}T_{k}(w_{2}-1)T_{l}(w_{3}-1)T_{m}(w_{1}-1)w_{1}^{k}w_{2}^{l}w_{3}^{m},\\
\label{m1}
&\sum_{k+l+m=n}\binom{n}{k,l,m}T_{k}(w_{3}-1)T_{l}(w_{1}-1)T_{m}(w_{2}-1)w_{2}^{k}w_{3}^{l}w_{1}^{m},\\
\label{n1}
&\sum_{k+l+m=n}\binom{n}{k,l,m}T_{k}(w_{3}-1)T_{l}(w_{2}-1)T_{m}(w_{1}-1)w_{1}^{k}w_{3}^{l}w_{2}^{m},\\
\label{o1}
&\sum_{k+l+m=n}\binom{n}{k,l,m}T_{k}(w_{2}-1)T_{l}(w_{1}-1)T_{m}(w_{3}-1)w_{3}^{k}w_{2}^{l}w_{1}^{m}.
\end{align}
\\
  However, (\ref{l1}) and (\ref{m1}) are equal to (\ref{f2}), as we
can see by applying the permutations $k\rightarrow l,l\rightarrow
m,m\rightarrow k$ for (\ref{l1}) and $k\rightarrow m,l\rightarrow
k,m\rightarrow l$ for (\ref{m1}). Similarly, we see that (\ref{n1})
and (\ref{o1}) are equal to (\ref{g2}), by applying permutations
$k\rightarrow l,l\rightarrow m,m\rightarrow k$ for (\ref{n1}) and
$k\rightarrow m,l\rightarrow k,m\rightarrow l$ for (\ref{o1}).
\begin{theorem}\label{A}
Let $w_{1},w_{2},w_{3}$ be any positive integers. Then the following
expression is invariant under any permutation of
$w_{1},w_{2},w_{3}$, so that it gives us six symmetries.
\begin{equation}\label{p1}
\begin{split}
&\sum_{k+l+m=n}\binom{n}{k,l,m}E_{k}(w_{1}y_{1})E_{l}(w_{2}y_{2})E_{m}(w_{3}y_{3})w_{1}^{l+m}w_{2}^{k+m}w_{3}^{k+l}\\
=&\sum_{k+l+m=n}\binom{n}{k,l,m}E_{k}(w_{1}y_{1})E_{l}(w_{3}y_{2})E_{m}(w_{2}y_{3})w_{1}^{l+m}w_{3}^{k+m}w_{2}^{k+l}\\
=&\sum_{k+l+m=n}\binom{n}{k,l,m}E_{k}(w_{2}y_{1})E_{l}(w_{1}y_{2})E_{m}(w_{3}y_{3})w_{2}^{l+m}w_{1}^{k+m}w_{3}^{k+l}\\
=&\sum_{k+l+m=n}\binom{n}{k,l,m}E_{k}(w_{2}y_{1})E_{l}(w_{3}y_{2})E_{m}(w_{1}y_{3})w_{2}^{l+m}w_{3}^{k+m}w_{1}^{k+l}\\
=&\sum_{k+l+m=n}\binom{n}{k,l,m}E_{k}(w_{3}y_{1})E_{l}(w_{1}y_{2})E_{m}(w_{2}y_{3})w_{3}^{l+m}w_{1}^{k+m}w_{2}^{k+l}\\
=&\sum_{k+l+m=n}\binom{n}{k,l,m}E_{k}(w_{3}y_{1})E_{l}(w_{2}y_{2})E_{m}(w_{1}y_{3})w_{3}^{l+m}w_{2}^{k+m}w_{1}^{k+l}.
\end{split}
\end{equation}
\end{theorem}
\begin{theorem}\label{B}
  Let $w_{1},w_{2},w_{3}$ be any odd positive integers. Then the
following expression is invariant under any permutation of
$w_{1},w_{2},w_{3}$, so that it gives us six symmetries.
\begin{equation}\label{q1}
\begin{split}
&\sum_{k+l+m=n}\binom{n}{k,l,m}E_{k}(w_{1}y_{1})E_{l}(w_{2}y_{2})T_{m}(w_{3}-1)w_{1}^{l+m}w_{2}^{k+m}w_{3}^{k+l}\\
=&\sum_{k+l+m=n}\binom{n}{k,l,m}E_{k}(w_{1}y_{1})E_{l}(w_{3}y_{2})T_{m}(w_{2}-1)w_{1}^{l+m}w_{3}^{k+m}w_{2}^{k+l}\\
=&\sum_{k+l+m=n}\binom{n}{k,l,m}E_{k}(w_{2}y_{1})E_{l}(w_{1}y_{2})T_{m}(w_{3}-1)w_{2}^{l+m}w_{1}^{k+m}w_{3}^{k+l}\\
=&\sum_{k+l+m=n}\binom{n}{k,l,m}E_{k}(w_{2}y_{1})E_{l}(w_{3}y_{2})T_{m}(w_{1}-1)w_{2}^{l+m}w_{3}^{k+m}w_{1}^{k+l}\\
=&\sum_{k+l+m=n}\binom{n}{k,l,m}E_{k}(w_{3}y_{1})E_{l}(w_{2}y_{2})T_{m}(w_{1}-1)w_{3}^{l+m}w_{2}^{k+m}w_{1}^{k+l}\\
=&\sum_{k+l+m=n}\binom{n}{k,l,m}E_{k}(w_{3}y_{1})E_{l}(w_{1}y_{2})T_{m}(w_{2}-1)w_{3}^{l+m}w_{1}^{k+m}w_{2}^{k+l}.
\end{split}
\end{equation}
\end{theorem}
Putting $w_{3}=1$ in (\ref{q1}), we get the following corollary.
\begin{corollary}\label{C}
Let $w_{1},w_{2}$ be any odd positive integers.
\begin{equation}\label{r1}
\begin{split}
&\sum_{k=0}^{n}\binom{n}{k}E_{k}(w_{1}y_{1})E_{n-k}(w_{2}y_{2})w_{1}^{n-k}w_{2}^{k}\\
=&\sum_{k=0}^{n}\binom{n}{k}E_{k}(w_{2}y_{1})E_{n-k}(w_{1}y_{2})w_{2}^{n-k}w_{1}^{k}\\
=&\sum_{k+l+m=n}\binom{n}{k,l,m}E_{k}(y_{1})E_{l}(w_{2}y_{2})T_{m}(w_{1}-1)w_{2}^{k+m}w_{1}^{k+l}\\
=&\sum_{k+l+m=n}\binom{n}{k,l,m}E_{k}(w_{2}y_{1})E_{l}(y_{2})T_{m}(w_{1}-1)w_{2}^{l+m}w_{1}^{k+l}\\
=&\sum_{k+l+m=n}\binom{n}{k,l,m}E_{k}(y_{1})E_{l}(w_{1}y_{2})T_{m}(w_{2}-1)w_{1}^{k+m}w_{2}^{k+l}\\
=&\sum_{k+l+m=n}\binom{n}{k,l,m}E_{k}(w_{1}y_{1})E_{l}(y_{2})T_{m}(w_{2}-1)w_{1}^{l+m}w_{2}^{k+l}.
\end{split}
\end{equation}
\end{corollary}
Letting further $w_{2}=1$ in (\ref{r1}), we have the following
corollary.
\begin{corollary}\label{D}
  Let $w_{1}$ be any odd positive integer.
\begin{equation}\label{s1}
\begin{split}
&\sum_{k=0}^{n}\binom{n}{k}E_{k}(w_{1}y_{1})E_{n-k}(y_{2})w_{1}^{n-k}\\
=&\sum_{k=0}^{n}\binom{n}{k}E_{k}(y_{1})E_{n-k}(w_{1}y_{2})w_{1}^{k}\\
=&\sum_{k+l+m=n}\binom{n}{k,l,m}E_{k}(y_{1})E_{l}(y_{2})T_{m}(w_{1}-1)w_{1}^{k+l}.
\end{split}
\end{equation}
\end{corollary}
\begin{theorem}\label{E}
Let $w_{1},w_{2},w_{3}$ be any odd positive integers. Then the
following expression  is invariant under any permutation of
$w_{1},w_{2},w_{3}$, so that it gives us six symmetries.
\begin{equation}\label{t1}
\begin{split}
&w_{1}^{n}\sum_{k=0}^{n}\binom{n}{k}E_{k}(w_{3}y_{1})\sum_{i=0}^{w_{1}-1}(-1)^iE_{n-k}(w_{2}y_{2}+\frac{w_{2}}{w_{1}}i)w_{3}^{n-k}w_{2}^{k}\\
=&w_{1}^{n}\sum_{k=0}^{n}\binom{n}{k}E_{k}(w_{2}y_{1})\sum_{i=0}^{w_{1}-1}(-1)^iE_{n-k}(w_{3}y_{2}+\frac{w_{3}}{w_{1}}i)w_{2}^{n-k}w_{3}^{k}\\
=&w_{2}^{n}\sum_{k=0}^{n}\binom{n}{k}E_{k}(w_{3}y_{1})\sum_{i=0}^{w_{2}-1}(-1)^iE_{n-k}(w_{1}y_{2}+\frac{w_{1}}{w_{2}}i)w_{3}^{n-k}w_{1}^{k}\\
=&w_{2}^{n}\sum_{k=0}^{n}\binom{n}{k}E_{k}(w_{1}y_{1})\sum_{i=0}^{w_{2}-1}(-1)^iE_{n-k}(w_{3}y_{2}+\frac{w_{3}}{w_{2}}i)w_{1}^{n-k}w_{3}^{k}\\
=&w_{3}^{n}\sum_{k=0}^{n}\binom{n}{k}E_{k}(w_{2}y_{1})\sum_{i=0}^{w_{3}-1}(-1)^iE_{n-k}(w_{1}y_{2}+\frac{w_{1}}{w_{3}}i)w_{2}^{n-k}w_{1}^{k}\\
=&w_{3}^{n}\sum_{k=0}^{n}\binom{n}{k}E_{k}(w_{1}y_{1})\sum_{i=0}^{w_{3}-1}(-1)^iE_{n-k}(w_{2}y_{2}+\frac{w_{2}}{w_{3}}i)w_{1}^{n-k}w_{2}^{k}.
\end{split}
\end{equation}
\end{theorem}
Letting $w_{3}=1$ in (\ref{t1}), we obtain alternative expressions
for the identities in (\ref{r1}).
\begin{corollary}\label{F}
  Let $w_{1},w_{2}$ be any odd positive integers.
\begin{equation}\label{u1}
\begin{split}
&\sum_{k=0}^{n}\binom{n}{k}E_{k}(w_{1}y_{1})E_{n-k}(w_{2}y_{2})w_{1}^{n-k}w_{2}^{k}\\
=&\sum_{k=0}^{n}\binom{n}{k}E_{k}(w_{2}y_{1})E_{n-k}(w_{1}y_{2})w_{2}^{n-k}w_{1}^{k}\\
=&w_{1}^{n}\sum_{k=0}^{n}\binom{n}{k}E_{k}(y_{1})\sum_{i=0}^{w_{1}-1}(-1)^iE_{n-k}(w_{2}y_{2}+\frac{w_{2}}{w_{1}}i)w_{2}^{k}\\
=&w_{1}^{n}\sum_{k=0}^{n}\binom{n}{k}E_{k}(w_{2}y_{1})\sum_{i=0}^{w_{1}-1}(-1)^iE_{n-k}(y_{2}+\frac{i}{w_{1}})w_{2}^{n-k}\\
=&w_{2}^{n}\sum_{k=0}^{n}\binom{n}{k}E_{k}(y_{1})\sum_{i=0}^{w_{2}-1}(-1)^iE_{n-k}(w_{1}y_{2}+\frac{w_{1}}{w_{2}}i)w_{1}^{k}\\
=&w_{2}^{n}\sum_{k=0}^{n}\binom{n}{k}E_{k}(w_{1}y_{1})\sum_{i=0}^{w_{2}-1}(-1)^iE_{n-k}(y_{2}+\frac{i}{w_{2}})w_{1}^{n-k}
\end{split}
\end{equation}
\end{corollary}
Putting further $w_{2}=1$ in (\ref{u1}), we have the alternative
expressions for the identities for (\ref{s1}).
\begin{corollary}\label{G}
  Let $w_{1}$ be any odd positive integer.
\begin{equation*}
\begin{split}
&\sum_{k=0}^{n}\binom{n}{k}E_{k}(y_{1})E_{n-k}(w_{1}y_{2})w_{1}^{k}\\
=&\sum_{k=0}^{n}\binom{n}{k}E_{k}(y_{2})E_{n-k}(w_{1}y_{1})w_{1}^{k}\\
=&w_{1}^{n}\sum_{k=0}^{n}\binom{n}{k}E_{k}(y_{1})\sum_{i=0}^{w_{1}-1}(-1)^{i}E_{n-k}(y_{2}+\frac{i}{w_{1}}).
\end{split}
\end{equation*}
\end{corollary}
\begin{theorem}\label{H}
Let $w_{1},w_{2},w_{3}$ be any odd positive integers. Then we have
the following three symmetries in $w_{1},w_{2},w_{3}$:
\begin{align}
\label{v1}
\sum_{k+l+m=n}\binom{n}{k,l,m}E_{k}(w_{1}y_{1})T_{l}(w_{2}-1)T_{m}(w_{3}-1)w_{1}^{l+m}w_{2}^{k+m}w_{3}^{k+l}\\
\label{w1}
=\sum_{k+l+m=n}\binom{n}{k,l,m}E_{k}(w_2y_{1})T_{l}(w_{3}-1)T_{m}(w_{1}-1)w_{2}^{l+m}w_{3}^{k+m}w_{1}^{k+l}\\
\label{x1}
=\sum_{k+l+m=n}\binom{n}{k,l,m}E_{k}(w_{3}y_{1})T_{l}(w_{1}-1)T_{m}(w_{2}-1)w_{3}^{l+m}w_{1}^{k+m}w_{2}^{k+l}.
\end{align}
\end{theorem}
Putting $w_{3}=1$ in (\ref{v1})-(\ref{x1}), we get the following
corollary.
\begin{corollary}\label{I}
  Let $w_{1},w_{2}$ be any odd positive integers.
\begin{equation}\label{y1}
\begin{split}
&\sum_{k=0}^{n}\binom{n}{k}E_{k}(w_{1}y_{1})T_{n-k}(w_{2}-1)w_{1}^{n-k}w_{2}^{k}\\
=&\sum_{k=0}^{n}\binom{n}{k}E_{k}(w_{2}y_{1})T_{n-k}(w_{1}-1)w_{2}^{n-k}w_{1}^{k}\\
=&\sum_{k+l+m=n}\binom{n}{k,l,m}E_{k}(y_{1})T_{l}(w_{1}-1)T_{m}(w_{2}-1)w_{1}^{k+m}w_{2}^{k+l}.
\end{split}
\end{equation}
\end{corollary}
Letting further $w_{2}=1$ in (\ref{y1}), we get the following
corollary. This is also obtained in [6, Cor.2] and mentioned in
\cite{D1}.
\begin{corollary}\label{J}
  Let $w_{1}$ be any odd positive integer.
\begin{equation}\label{z1}
E_{n}(w_{1}y_{1})=\sum_{k=0}^{n}\binom{n}{k}E_{k}(y_{1})T_{n-k}(w_{1}-1)w_{1}^{k}.
\end{equation}
\end{corollary}
\begin{theorem}\label{K}
Let $w_{1},w_{2},w_{3}$ be any odd positive integers. Then the
following expression  is invariant under any permutation of
$w_{1},w_{2},w_{3}$, so that it gives us six symmetries.
\begin{equation}\label{a2}
\begin{split}
&w_{1}^{n}\sum_{k=0}^{n}\binom{n}{k}\sum_{i=0}^{w_{1}-1}(-1)^iE_{k}(w_{2}y_{1}+\frac{w_{2}}{w_{1}}i)T_{n-k}(w_{3}-1)w_{2}^{n-k}w_{3}^{k}\\
=&w_{1}^{n}\sum_{k=0}^{n}\binom{n}{k}\sum_{i=0}^{w_{1}-1}(-1)^iE_{k}(w_{3}y_{1}+\frac{w_{3}}{w_{1}}i)T_{n-k}(w_{2}-1)w_{3}^{n-k}w_{2}^{k}\\
=&w_{2}^{n}\sum_{k=0}^{n}\binom{n}{k}\sum_{i=0}^{w_{2}-1}(-1)^iE_{k}(w_{1}y_{1}+\frac{w_{1}}{w_{2}}i)T_{n-k}(w_{3}-1)w_{1}^{n-k}w_{3}^{k}\\
=&w_{2}^{n}\sum_{k=0}^{n}\binom{n}{k}\sum_{i=0}^{w_{2}-1}(-1)^iE_{k}(w_{3}y_{1}+\frac{w_{3}}{w_{2}}i)T_{n-k}(w_{1}-1)w_{3}^{n-k}w_{1}^{k}\\
=&w_{3}^{n}\sum_{k=0}^{n}\binom{n}{k}\sum_{i=0}^{w_{3}-1}(-1)^iE_{k}(w_{1}y_{1}+\frac{w_{1}}{w_{3}}i)T_{n-k}(w_{2}-1)w_{1}^{n-k}w_{2}^{k}\\
=&w_{3}^{n}\sum_{k=0}^{n}\binom{n}{k}\sum_{i=0}^{w_{3}-1}(-1)^iE_{k}(w_{2}y_{1}+\frac{w_{2}}{w_{3}}i)T_{n-k}(w_{1}-1)w_{2}^{n-k}w_{1}^{k}.
\end{split}
\end{equation}
\end{theorem}
Putting $w_{3}=1$ in (\ref{a2}), we obtain the following corollary.
In Section 1, the identities in (\ref{y1}), (\ref{b2}), and
(\ref{d2}) are combined to give those in (\ref{h})-(\ref{o}).
\begin{corollary}\label{L}
  Let $w_{1},w_{2}$ be any odd positive integers.
\begin{equation}\label{b2}
\begin{split}
&w_{1}^{n}\sum_{i=0}^{w_{1}-1}(-1)^{i}E_{n}(w_{2}y_{1}+\frac{w_{2}}{w_{1}}i)\\
=&w_{2}^{n}\sum_{i=0}^{w_{2}-1}(-1)^{i}E_{n}(w_{1}y_{1}+\frac{w_{1}}{w_{2}}i)\\
=&\sum_{k=0}^{n}\binom{n}{k}E_{k}(w_{2}y_{1})T_{n-k}(w_{1}-1)w_{2}^{n-k}w_{1}^{k}\\
=&\sum_{k=0}^{n}\binom{n}{k}E_{k}(w_{1}y_{1})T_{n-k}(w_{2}-1)w_{1}^{n-k}w_{2}^{k}\\
=&w_{1}^{n}\sum_{k=0}^{n}\binom{n}{k}\sum_{i=0}^{w_{1}-1}(-1)^{i}E_{k}(y_{1}+\frac{i}{w_{1}})T_{n-k}(w_{2}-1)w_{2}^{k}\\
=&w_{2}^{n}\sum_{k=0}^{n}\binom{n}{k}\sum_{i=0}^{w_{2}-1}(-1)^{i}E_{k}(y_{1}+\frac{i}{w_{2}})T_{n-k}(w_{1}-1)w_{1}^{k}.
\end{split}
\end{equation}
\end{corollary}
Letting further $w_{2}=1$ in (\ref{b2}), we get the following
corollary. This is the multiplication formula for Euler polynomials
together with the relatively new identity mentioned in (\ref{z1}).
\begin{corollary}\label{M}
  Let $w_{1}$ be any odd positive integer.
\begin{equation*}
\begin{split}
E_{n}(w_{1}y_{1})&=w_{1}^{n}\sum_{i=0}^{w_{1}-1}(-1)^{i}E_{n}(y_{1}+\frac{i}{w_{1}})\\
&=\sum_{k=0}^{n}\binom{n}{k}E_{k}(y_{1})T_{n-k}(w_{1}-1)w_{1}^{k}.
\end{split}
\end{equation*}
\end{corollary}
\begin{theorem}\label{N}
Let $w_{1},w_{2},w_{3}$ be any odd positive integers. Then we have
the following three symmetries in $w_{1},w_{2},w_{3}$:
\begin{equation}\label{c2}
\begin{split}
&(w_{1}w_{2})^{n}\sum_{i=0}^{w_{1}-1}\sum_{j=0}^{w_{2}-1}(-1)^{i+j}E_{n}(w_{3}y_{1}+\frac{w_{3}}{w_{1}}i+\frac{w_{3}}{w_{2}}j)\\
=&(w_{2}w_{3})^{n}\sum_{i=0}^{w_{2}-1}\sum_{j=0}^{w_{3}-1}(-1)^{i+j}E_{n}(w_{1}y_{1}+\frac{w_{1}}{w_{2}}i+\frac{w_{1}}{w_{3}}j)\\
=&(w_{3}w_{1})^{n}\sum_{i=0}^{w_{3}-1}\sum_{j=0}^{w_{1}-1}(-1)^{i+j}E_{n}(w_{2}y_{1}+\frac{w_{2}}{w_{3}}i+\frac{w_{2}}{w_{1}}j).
\end{split}
\end{equation}
\end{theorem}
Letting $w_{3}=1$ in (\ref{c2}), we have the following corollary.
\begin{corollary}\label{O}
  Let $w_{1},w_{2}$ be any odd positive integers.
\begin{equation}\label{d2}
\begin{split}
&w_{1}^{n}\sum_{j=0}^{w_{1}-1}(-1)^{j}E_{n}(w_{2}y_{1}+\frac{w_{2}}{w_{1}}j)\\
=&w_{2}^{n}\sum_{i=0}^{w_{2}-1}(-1)^{i}E_{n}(w_{1}y_{1}+\frac{w_{1}}{w_{2}}i)\\
=&(w_{1}w_{2})^{n}\sum_{i=0}^{w_{1}-1}\sum_{j=0}^{w_{2}-1}(-1)^{i+j}E_{n}(y_{1}+\frac{i}{w_{1}}+\frac{j}{w_{2}}).
\end{split}
\end{equation}
\end{corollary}
\begin{theorem}\label{P}
Let  $w_{1},w_{2},w_{3}$ be any positive integers. Then we have the
following two symmetries in  $w_{1},w_{2},w_{3}$:
\begin{equation}\label{e2}
\begin{split}
&\sum_{k+l+m=n}\binom{n}{k,l,m}E_{k}(w_{1}y)E_{l}(w_{2}y)E_{m}(w_{3}y)w_{3}^{k}w_{1}^{l}w_{2}^{m}\\
=&\sum_{k+l+m=n}\binom{n}{k,l,m}E_{k}(w_{1}y)E_{l}(w_{3}y)E_{m}(w_{2}y)w_{2}^{k}w_{1}^{l}w_{3}^{m}.
\end{split}
\end{equation}
\end{theorem}
\begin{theorem}\label{Q}
Let $w_{1},w_{2},w_{3}$ be any odd positive integers. Then we have
the following two symmetries in $w_{1},w_{2},w_{3}$:
\begin{align}
\label{f2}
&\sum_{k+l+m=n}\binom{n}{k,l,m}T_{k}(w_{1}-1)T_{l}(w_{2}-1)T_{m}(w_{3}-1)w_{3}^{k}w_{1}^{l}w_{2}^{m}\\
\label{g2}
=&\sum_{k+l+m=n}\binom{n}{k,l,m}T_{k}(w_{1}-1)T_{l}(w_{3}-1)T_{m}(w_{2}-1)w_{2}^{k}w_{1}^{l}w_{3}^{m}.
\end{align}
\\
\end{theorem}
Putting $w_{3}=1$ in (\ref{f2}) and (\ref{g2}),  we get the
following corollary.
\begin{corollary}\label{R}
Let $w_{1},w_{2}$ be any odd positive integers.
\begin{equation*}
\begin{split}
&\sum_{k=0}^{n}\binom{n}{k}T_{k}(w_{2}-1)T_{n-k}(w_{1}-1)w_{1}^{k}\\
=&\sum_{k=0}^{n}\binom{n}{k}T_{k}(w_{1}-1)T_{n-k}(w_{2}-1)w_{2}^{k}.
\end{split}
\end{equation*}
\end{corollary}
\bibliographystyle{amsplain}

\end{document}